\documentclass[12pt]{article}
\usepackage{amsfonts,amsmath}
\def \beq {\begin{eqnarray}}
\def \eeq {\end{eqnarray}}
\def \beqn {\begin{eqnarray*}}
\def \eeqn {\end{eqnarray*}}

\newcommand{\halmos}{\rule{1ex}{1.4ex}}

\newcounter{for}[section]

\newtheorem{itlemma}{Lemma}[section]
\newtheorem{itproposition}[itlemma]{Proposition}
\newtheorem{theorem}[itlemma]{Theorem}
\newtheorem{itcorollary}[itlemma]{Corollary}
\newtheorem{itremark}[itlemma]{Remark}
\newtheorem{itremarks}[itlemma]{Remarks}
\newtheorem{itdefinition}[itlemma]{Definition}
\newtheorem{itexample}[itlemma]{Example}

\newenvironment{fact}{\begin{itfact}\rm}{\end{itfact}}
\newenvironment{claim}{\begin{itclaim}\rm}{\end{itclaim}}
\newenvironment{lemma}{\begin{itlemma}}{\end{itlemma}}
\newenvironment{remark}{\begin{itremark}\rm}{\end{itremark}}
\newenvironment{remarks}{\begin{itremarks} \rm}{\end{itremarks}}
\newenvironment{corollary}{\begin{itcorollary}}{\end{itcorollary}}
\newenvironment{proposition}{\begin{itproposition}}{\end{itproposition}}
\newenvironment{definition}{\begin{itdefinition}\rm}{\end{itdefinition}}
\newenvironment{example}{\begin{itexample}\rm}{\end{itexample}}
\newenvironment{proof}{\noindent {\em Proof}.\ \
}{\hspace*{\fill}$\halmos$\medskip}
\newcommand{\be}[1]{\addtocounter{for}{1} \begin{equation}\label{#1}}
\newcommand{\ee}{\end{equation}}
\newcommand{\bl}[1]{\begin{lemma}\label{#1}}
\newcommand{\br}[1]{\begin{remark}\label{#1}}
\newcommand{\brs}[1]{\begin{remarks}\label{#1}}
\newcommand{\bt}[1]{\begin{theorem}\label{#1}}
\newcommand{\bd}[1]{\begin{definition}\label{#1}}
\newcommand{\bp}[1]{\begin{proposition}\label{#1}}
\newcommand{\bc}[1]{\begin{corollary}\label{#1}}
\newcommand{\bfact}[1]{\begin{fact}\label{#1}}
\newcommand{\bex}[1]{\begin{example}\label{#1}}
\newcommand{\ec}{\end{corollary}}
\newcommand{\efact}{\end{fact}}
\newcommand{\eex}{\end{example}}
\newcommand{\el}{\end{lemma}}
\newcommand{\er}{\end{remark}}
\newcommand{\ers}{\end{remarks}}
\newcommand{\et}{\end{theorem}}
\newcommand{\ed}{\end{definition}}
\newcommand{\ep}{\end{proposition}}
\newcommand{\epr}{\end{proof}}
\newcommand{\bpr}{\begin{proof}}
\newcommand{\bcl}[1]{\begin{claim}\label{#1}}
\newcommand{\ecl}{\end{claim}}

\newcommand{\ecs}{\end{corollary}}
\newcommand{\eers}{\end{exercise}}
\newcommand{\eexs}{\end{example}}
\newcommand{\eems}{\end{example}}
\newcommand{\els}{\end{lemma}}
\newcommand{\eles}{\end{lemmaex}}
\newcommand{\ets}{\end{theorem}}
\newcommand{\eds}{\end{definition}}
\newcommand{\eps}{\end{proposition}}

\newcommand{\bi}{\begin{itemize}}
\newcommand{\ei}{\end{itemize}}
\newcommand{\ben}{\begin{enumerate}}
\newcommand{\een}{\end{enumerate}}

\def\vbar{\mathchoice{\vrule height6.3ptdepth-.5ptwidth.8pt\kern-.8pt}
   {\vrule height6.3ptdepth-.5ptwidth.8pt\kern-.8pt}
   {\vrule height4.1ptdepth-.35ptwidth.6pt\kern-.6pt}
   {\vrule height3.1ptdepth-.25ptwidth.5pt\kern-.5pt}}
\def\fudge{\mathchoice{}{}{\mkern.5mu}{\mkern.8mu}}
\def\bbc#1#2{{\rm \mkern#2mu\vbar\mkern-#2mu#1}}
\def\bbb#1{{\rm I\mkern-3.5mu #1}}
\def\bba#1#2{{\rm #1\mkern-#2mu\fudge #1}}
\def\bb#1{{\count4=`#1 \advance\count4by-64 \ifcase\count4\or\bba A{11.5}\or
   \bbb B\or\bbc C{5}\or\bbb D\or\bbb E\or\bbb F \or\bbc G{5}\or\bbb H\or
   \bbb I\or\bbc J{3}\or\bbb K\or\bbb L \or\bbb M\or\bbb N\or\bbc O{5} \or
   \bbb P\or\bbc Q{5}\or\bbb R\or\bbc S{4.2}\or\bba T{10.5}\or\bbc U{5}\or
   \bba V{12}\or\bba W{16.5}\or\bba X{11}\or\bba Y{11.7}\or\bba Z{7.5}\fi}}

\def \da{\downarrow}
\def \ua{\uparrow}

\def \P{{\cal{P}}}

\def \A {{\cal{A}}}

\def \FF {{\cal{F}}}

\def \E {{\cal{E}}}
\def \C {{\cal{C}}}
\def \D {{\cal{D}}}
\def \Ra {{\cal{R}}}

\newcommand{\ba}[1]{\addtocounter{for}{1} \begin{eqnarray}\label{#1}}
\newcommand{\ea}{\end{eqnarray}}

\def\sqr#1#2{{\vcenter{\vbox{\hrule height .#2pt
                             \hbox{\vrule width .#2pt height#1pt \kern#1pt
                                   \vrule width .#2pt}
                             \hrule height .#2pt}}}}

\def\pmb#1{\setbox0=\hbox{#1}%
   \kern-.025em\copy0\kern-\wd0
   \kern.05em\copy0\kern-\wd0
   \kern-.025em\raise.0433em\box0 }
\def\sqr#1#2{{\vcenter{\vbox{\hrule height.#2pt
     \hbox{\vrule width.#2pt height#1pt \kern#1pt
   \vrule width.#2pt}\hrule height.#2pt}}}}

\def\reff#1{(\ref{#1})}

\def\N{{\mathbb N}}   
\def\Z{{\mathbb Z}}
\def\R{{\mathbb R}}

\def\P{{\mathbb P}}
\def\E{{\mathbb E}}

\def\e{\epsilon}                

\def\e{\epsilon}
\def\d{\delta}
\def\l{\lambda}

\def\a{\alpha}

\def\v{\varphi}

\newcommand {\pare}[1] {\left( {#1} \right)}
\newcommand {\cro}[1] {\left[ {#1} \right]}
\newcommand {\nor}[1] { \left\| {#1} \right\|}
\newcommand {\refeq}[1] {(\ref{#1})}
\def \ind {\hbox{ 1\hskip -3pt I}}
\newcommand {\acc}[1] {\left\{ {#1} \right\}}

\begin{document}

\title{A Note on Random Walk in Random Scenery.}
\author{Amine Asselah \& Fabienne Castell\\C.M.I., Universit\'e de Provence,\\
39 Rue Joliot-Curie, \\F-13453 Marseille cedex 13, France\\
asselah@cmi.univ-mrs.fr \& castell@cmi.univ-mrs.fr}
\date{}
\maketitle
\begin{abstract}
We consider a random walk in random scenery $\{X_n=
\eta(S_0)+\dots+\eta(S_n),n\in \N\}$,
where a centered walk $\{S_n,n\in \N\}$ is independent of the scenery 
$\{\eta(x),x\in \Z^d\}$, consisting of
symmetric i.i.d.\ with tail distribution $P(\eta(x)>t)\sim \exp(-c_{\a}t^{\a})$,
with $1\le \a<d/2$. We study the probability, when averaged over both 
randomness, that $\{X_n>n y\}$ for $y>0$, and $n$ large. 
In this note, 
we show that the large deviation estimate is of order $\exp(-c (ny)^a)$,
with $a=\a/(\a+1)$.
\end{abstract}



{\em Running head}: Random Walk in Random Scenery.

\section{Introduction}
We consider a centered random walk $\{S_k,k\in \N\}$ on $\Z^d$. When
$S_0=x$, we denote the law of the walk by
$\P_{x}$ and the expectation with respect to this law by $\E_x$.
Each site $x\in \Z^d$ is associated with
a random variable $\eta(x)$, and we assume that the {\it scenery}
$\{\eta(x),x\in \Z^d\}$ consists of symmetric
i.i.d.\ unbounded random variables, independent of the random walk.
We denote the law of the scenery by $P_{\eta}$, and by $E_{\eta}$ the
expectation with respect to this law.

The random walk in random scenery (RWRS) is the process
$\{X_n, n \in \N\}$ defined by
\be{rwrs.def}
X_n := \sum_{k=0}^{n} \eta(S_k)=\sum_{x\in \Z^d} l_n(x)\eta(x),\quad
\text{where}\quad l_n(x)=\sum_{k=0}^{n}\ind_{S_k=x} \, .
\ee
RWRS has been introduced by Kesten, Spitzer \cite{KS}, and Borodin
\cite{B1,B2} as a case-study for sums
of dependent random variables, in order to exhibit new scaling and new
self-similar limiting laws. Indeed, the convergence in law of $X_n$,
studied for $d \neq 2$ in \cite{KS,B1,B2}, and for
$d=2$ by Bolthausen~\cite{B}, needs a {\it super-diffusive} scaling in
dimensions 1 and 2. In terms of the mean square, the dominant orders 
are the following.
\be{TLC}
\E_0 \otimes E_{\eta} \cro{ X_n ^2}
\simeq \left\{ \begin{array}{ll}
                n^{3/2} & \mbox{ for } d=1 \, ,
                \\
                n \log(n)  & \mbox{ for } d=2 \, ,
                \\
                n  & \mbox{ for } d \geq 3 \, .
                \end{array}
        \right.
\ee

Recently, the moderate and large deviations for $X_n$ have been
studied in~\cite{CP,AC1,AC2,C} for the Brownian motion
in various sceneries, and in~\cite{GKS,HGK} in the original random walk setting.

We distinguish three regimes depending on the tail parameter $\alpha$ of
the scenery variable
\be{eq1.1}
P_{\eta}(\eta(x)>t)\sim \exp(-c_{\a} t^{\a}),\quad\text{for $t$ large}.
\ee
\begin{itemize}
\item
When $\a<1$, $\eta(x)$ has no exponential moments, and is called
a heavy-tail variable. In a recent paper~\cite{HGK},
van der Hofstad, Gantert and K\"onig deal with $X_n=\sum l_n(x)\eta(x)$, by
conditioning on the local times $\{l_n(x),x\in \Z^d\}$, thus obtaining a
weighted sum of i.i.d.\ heavy-tail variables. Using classical heavy-tail
estimates (as those of \cite{Nagaev}), they show that $\{X_n>ny\}$
is realized as only one term of the series
reaches level $ny$. Thus, in terms of logarithmic equivalence ($\approx$),
\[
\P_0 \otimes P_{\eta} \cro{ X_n  \geq n y}
\approx \P_0 \otimes P_{\eta} \cro{l_n(0) \eta(0) \geq n y}.
\]
Now, recall that for a time $k\ll n$, 
the local time at site 0, $l_n(0)$, satisfies the following property
\be{eq1.6}
\P_0(l_n(0)=k)\approx \exp(-\kappa_0 k),\quad\text{with}\quad
\kappa_0:=\log\left(\frac{1}{\P_0(H_0<\infty)}\right) \, ,
\ee
where $H_0 = \inf\acc{n \geq 1, X_n =0}$.
Thus, for $y>0$, \cite{HGK} shows that for an explicit $J>0$
\ba{eq1.7}
\P_0\otimes P_{\eta}\left(l_n(0)\eta(0)>ny\right)
&\approx &\sup_{k=1,2,\dots} \P_0\left(l_n(0)=k\right)
P_{\eta}\left(\eta(0)>\frac{ny}{k}\right)\cr
&\approx & \exp\left(- \inf_{k\ge 1}\left(\kappa k+c_{\a}
(\frac{ny}{k})^{\a}\right)\right)\cr
&\approx& \exp\left(-J (ny)^{\frac{\a}{\a+1}}\right).
\ea
In the optimal strategy $l_n(0)$ is of order $(ny)^{\a/(\a+1)}$.
\item When $\alpha > \max(d/2,1)$, a different behavior holds: for all $y>0$
\be{GDlight}
\P_0 \otimes P_{\eta}
\cro{ X_n  \geq n y} \approx\exp \pare{- n^{\frac{d}{d+2}} J(y)},\quad
\text{( with a $J(y)>0$ known explicitly)}. 
\ee
This result is proved
in \cite{AC2} for Brownian motion in a bounded scenery (i.e.\ $\alpha = \infty$), in
\cite{C} for a Gaussian scenery ($\alpha =2$ and $d\le 3$), and in \cite{GKS}
for a random walk in a general scenery. 
The best strategy to realize $\{X_n>ny\}$ is the following.
\begin{itemize}
\item Force the random walk to spend all its time in a ball of radius $r_n$
with $1\ll r_n\ll \sqrt{n}$, in such a way that for $x$ in this ball, $l_n(x)$
is of order $n/r_n^d$. 
This has a cost of order $\exp(- n/r_n^2)$
\item Require the scenery to satisfy $
\sum_{\nor{x} \leq r_n} \eta(x) \geq  r_n^dy$. This has a cost of order
$\exp(-r_n^d)$.
\end{itemize}
The exponent $d/(d+2)$ appears as one sets equal $n/r_n^2$ and $r_n^d$.
Thus, in the optimal strategy, the walk spends a time $n^{2/(d+2)}$
on each site of a ball of about $n^{d/(d+2)}$ sites. 

Observe that when $\eta(x)$ satisfies \reff{eq1.1},
then $l_n(x)\eta(x)$ has a heavy tail (see \reff{eq1.7}). Also
$X_n=\sum \eta(x) l_n(x)$ is a sum of about $n$-terms 
(in dimensions $d\ge 3$).
However, the variables $\{l_n(x)\eta(x),x \in \Z^d\} \}$
are not independent, and the extreme value of the sum does not dominate.

\item The regime $1\le \a< d/2$ is the purpose of this note.
\end{itemize}

Our main result is the following. 
\bp{prop1.1} Let $\{S_n,n\in \N\}$ be a walk with centered independent increments
with finite exponential moments. Assume that $\{\eta(x),x\in \Z^d\}$ are symmetric
i.i.d.\ variables with tail parameter $\a$ with $1\le \a< d/2$, and whose
law has a density decreasing on $\R^+$. 
There are $c_1,c_2>0$, such that when $n$ is large enough
\be{LD.ineq}
\exp(- c_1 (ny)^{\a/(\a+1)}) \leq
\log \P_0\otimes P_{\eta} \cro{X_n \geq  n y}\leq\exp(- c_2 (ny)^{\a/(\a+1)}).
\ee
\ep
In the course of deriving the upper bound, we rely on a localization lemma
of independent interest. 

\bl{lem-amin} Assume $d\ge 3$. There is a constant $\kappa_d>0$ such
that for any $\Lambda\subset \Z^d$, and any $t>0$
\be{key-amin}
\P_0(l_{\infty}(\Lambda)>t)\le
\exp\left(-\kappa_d \frac{t}{|\Lambda|^{2/d}}\right),
\ee
where $l_{\infty}(\Lambda)$ is the total sojourn
time of the walk in the region $\Lambda$.
\el
Note that the recent paper \cite{CFRRS}
gives a representation of $\P_0(l_{\infty}(\Lambda)>t)$, in terms
of the eigenvalues and eigenvectors of the matrix whose entries
are the Green function restricted to $\Lambda$. It is not 
clear to us how to deduce Lemma~\ref{lem-amin} from the type of
representation of \cite{CFRRS}.

This note is organized as follows.
We specify the model in Section~\ref{model}.
In Section~\ref{lower}, we deal with the lower bound.
In Section~\ref{upper}, we deal with the upper bound. 
Finally, we have gathered in the Appendix
the proof of Lemma~\ref{lem-amin}, and the proof of some technical facts.

\section{Model}
\label{model}

\noindent
{\bf  Assumptions on the random walk}. We assume that the increments of
the walk are centered, with finite exponential moments, i.e.
\be{RW.hyp}
S_k = \sum_{j=1}^{k} \xi_j\, ,  \,\,
\xi_j \mbox{ i.i.d  }  \, , \,\,
E\cro{\xi_1} = 0 \, , \,\,
E\cro{\exp(\lambda \xi_1)} < \infty \mbox{ for all } \lambda \in \R^d \, .
\ee
It is then easy to see that there exist constants $C,c > 0$, such that
for all $n$,
\be{ineq-exit}
\P_0\cro{\max_{k\leq n} \nor{S_k} \geq n} \leq C \exp(-cn)  \, .
\ee
\noindent
{\bf  Assumptions on the scenery}.
Besides our basic tail assumption \refeq{eq1.1},
we make assumptions on the law of the scenery whose goal is
to simplify the technical parts. Thus, we say that a random variable
with value in $\R$ is {\it bell-shaped}, if its law has
a density with respect to Lebesgue which is even, and decreasing on $\R^+$.
Throughout the paper, we will assume that
$\{\eta(x); x \in \Z^d\}$ are i.i.d\ and bell-shaped,
with following handy consequence, proved in the Appendix.
\bl{lem-handy}
When $\{\eta(x),x\in \Z^d\}$ have independent bell-shaped densities,
then for any $\Lambda$ finite subset of $\Z^d$, and any $y>0$
\be{ineq-handy}
P\left(\sum_{x\in\Lambda} \a_x \eta(x)>y\right)\le
P\left(\sum_{x\in\Lambda} \beta_x \eta(x)>y\right),
\quad\text{if}\quad 0\le \a_x\le \beta_x\text{ for all }x\in \Lambda.
\ee
\el
A typical use of Lemma~\ref{lem-handy} is the following bound
\be{ineq-use}
P\left(\sum_{\Lambda} \eta(x)>\frac{y}{\min \a_x}\right)\le
P\left(\sum_{\Lambda} \a_x\eta(x)>y\right)\le
P\left(\sum_{\Lambda} \eta(x)>\frac{y}{\max \a_x}\right).
\ee

\noindent {\bf Some notations.}
Throughout the paper, we set
$a := \alpha/(\alpha+1)$ and $b := 1/(\alpha + 1)$,
and for $x \in \Z^d$, $\nor{x} := \max_{i =1,\cdots, d} |x_i|$.
Finally, when considering the variables $\{\eta(x), x \in \Lambda\}$ for
a finite region $\Lambda$ of cardinality $L$, we will sometimes use
the notation $\{\eta_j, 1 \leq j \leq L\}$.

\section{Lower Bound}\label{lower}
We show in this section the following simple estimate.
\bl{lem-lower}
There is a constant $c_1>0$ such that, for any $y>0$, and $n$ large
\be{eqlb.4}
P\left( \sum_{x\in \Z^d} l_n(x) \eta(x)>ny\right)\ge \exp\left(-c_1 (ny)^a
\right).
\ee
\el
\bpr
The bound \reff{eqlb.4} is obtained by using Lemma~\ref{lem-handy}. Thus,
\ba{eqlb.1}
P\left( \sum_x l_n(x) \eta(x)>ny\right)&\ge&
P\left( l_n(0) \eta(0)>ny\right)=\sum_{k>0}
\P_0(l_n(0)=k)P_{\eta}\left(\eta(0)>\frac{ny}{k}\right)\cr
&\ge& \P_0(l_n(0)=k)P_{\eta}\left(\eta(0)>\frac{ny}{k}\right)\quad
\text{for any }k.
\ea
We choose an $n$-depending $k$, for instance $k_n=[(ny)^a]$. 
Since $k_n/(ny)^{a}\to 1$ as $n$ tends to infinity, we have for any
$\e>0$ and $n$ large that
\beqn
P_{\eta}\left(\eta(0)>\frac{ny}{k_n}\right)&\ge &\exp
\left( -c_{\a}(1+\e)\left(\frac{ny}{k_n}\right)^{\a}\right)\cr
&\ge &\exp \left( -c_{\a}(1+2\e)(ny)^a\right).
\eeqn
Now, if we set $\kappa_0=\log(1/\P_0(H_0<\infty))$, then
\ba{eqlb.2}
\P_0(l_n(0)=k_n)&\ge&\P_0(H_0\le \frac{n}{k_n})^{k_n}\cr
&\ge& \left( e^{-\kappa_0}-\P_0(\frac{n}{k_n}<H_0<\infty)\right)^{k_n}.
\ea
Thus, for any $\kappa>\kappa_0$, we have for $n$ large
enough
\be{eqlb.3}
P\left( \sum_x l_n(x) \eta(x)>ny\right)\ge \exp\left(-
(\kappa+(1+2\e)c_{\a}) (ny)^{a}\right)
\ee
This concludes the proof.
\epr
\section{Upper Bound}\label{upper}
The case $\a=1$ is special and much simpler than $\a>1$. Thus, we will treat
the former specifically in Remark~\ref{rem-a1}. Henceforth,
we assume that $\a>1$ and we recall that $a:=\a/(\a+1)$, and $b=1-a<a$. We consider
a subdivision of $[b,a]$, $b_1=b<b_2<\dots<b_{N+1}=a$,
and a decomposition of $y>0$ into 
positive constants $\{y_{\da},y_0,\dots,y_N,y_{\ua}\}$ summing
up to $y$. We will specify $N$ and $\{b_i,y_i,i=1,\dots,N\}$ 
after we partition the range of the walk $\Ra_n$, into $N+3$ sets.
For $1\le i\le N$, we set
\be{equb.1}
\D_i:=\{x\in \Ra_n: y^an^{b_i}\le l_n(x)< y^an^{b_{i+1}}\},
\ee
and for a {\it small} constant $z$ to be chosen later
\be{def-D0}
\D_0:=\{x\in \Ra_n: zn^{b}\le l_n(x)< y^an^{b}\},
\ee
and lastly, for the two sets at the extremities
\be{def-DUL}
\D_{\da}:=\{x\in \Ra_n: l_n(x)< zn^b\},\quad\text{and}\quad
\D_{\ua}:=\{x\in \Ra_n: l_n(x)\ge (yn)^a\}.
\ee
Thus,
\be{equb.2}
\{\sum_{x\in \Z^d} l_n(x) \eta(x)>n y\}\subset \bigcup_{i=0}^N \{ \sum_{x\in\D_i}
l_n(x)\eta(x)> n y_i\}\bigcup\{ \sum_{x\in\D_{\da}}
l_n(x)\eta(x)> n y_{\da}\}\bigcup\{ \D_{\ua}\not=\emptyset\}.
\ee
Thus, if we define $\A:=\{\max_{k\le n}||S_k||<n\}$, and recall that
$\P_0(\A^c)$ is negligible compared to $\exp(-n^a)$ by \reff{ineq-exit}, then
\ba{equb.3}
P\left(\sum_{x\in \Z^d} l_n(x) \eta(x)>n y\right)&\le&\P_0(\A^c)+
\sum_{i=0}^N P\left(\A, \sum_{x\in\D_i} l_n(x) \eta(x)> n y_i\right)\cr
&&+P\left(\A, \sum_{x\in\D_{\da}} l_n(x) \eta(x)> n y_{\da}\right)+
\P_0\left(\D_{\ua}\not=\emptyset\right).
\ea
We will now estimate each terms separately in the next Section.
However, in the course of obtaining an upper bound, we will fall
on the following requirement: 
we will need a positive $\beta$, independent of $n$,
such that for $i\le N$
\be{ineq.req1}
\beta y\le y_i n^{(a-b_{i+1})-(1-\d_0)(a-b_i)}\quad
\text{with the positive constant}\quad \d_0:=\frac{1/\a-2/d}{1-2/d}<1.
\ee
Thus, a simple choice of $\{b_i,y_i\}$ which fulfills \reff{ineq.req1}
is $y_N:=\beta y n^{(1-\delta_0)(a-b_N)}$, and
\be{def-epsilon}
\forall i<N,\quad y_i=\beta yn^{-\e_0(a-b_{i+1})},\text{  and  }
(a-b_i)=(1+\e_0)(a-b_{i+1}),\quad\text{with}\quad
\e_0=\frac{\d_0/2}{1-\d_0/2}.
\ee
To explicit further the choices in \reff{def-epsilon}, we introduce more notations: 
\be{def-z}
z_1=a-b_N,\quad z_2=b_N-b_{N-1},\dots,z_N=b_2-b_1.
\ee
Thus, \reff{def-epsilon} is fulfilled when $z_2=\e_0 z_1$ and for $i>2$
\be{choice}
z_i=\e_0(z_1+\dots+z_{i-1})=(1+\e_0)z_{i-1}=(1+\e_0)^{i-2}z_2=
(1+\e_0)^{i-2}\e_0 z_1.
\ee
Note that for $i< N$
\[
a-b_{i+1}=z_1+\dots+z_{N-i}=\frac{z_{N-i+1}}{\e_0}=(1+\e_0)^{N-i-1}z_1.
\]
The condition on $y_i$ in \reff{def-epsilon} will be fulfilled if we choose
$z_1=\chi /\log(n)$, for a constant $\chi$ to be tuned later. 
Indeed, we obtain $y_N=\beta y \exp(\chi(1-\delta_0))$, and 
\be{choice-y}
\forall i<N,\quad y_i=\beta y 
\exp\left( -\log(n)z_1 \e_0(1+\e_0)^{N-i-1}\right)=
\beta y \exp\left( -\chi \e_0(1+\e_0)^{N-i-1}\right).
\ee
Thus, since
\be{need-choice}
\sum_{i=1}^{\infty} (1+\e_0)^i=\infty,\quad\text{and}\quad
\sum_{i=1}^{\infty} \exp\left( -\chi \e_0(1+\e_0)^{i}\right)<\infty,
\ee
one can find $N$ finite,
of order $\log(\log(n))$, 
$\chi>0$ and $\beta>0$ independent of $n$ (or rather $\chi_n$ and $\beta_n$
can be chosen to converge to positive constants, and we omit the subscript $n$)
such that
\be{good-choice}
a-b=\sum_{i=1}^{N} z_i=\frac{\chi}{\log(n)}(1+\e_0)^{N-1},\quad\text{and}\quad
\sum_{i=0}^{N}y_i+y_{\da}+y_{\ua}=y.
\ee
Actually, the choice of $y_{\da},y_{\ua}$ is arbitrary since we are
not after the exact constant in front of the speed $n^a$.
For instance, we choose $y_{\da}=y_{\ua}=y/3$.
\subsection{Contribution of $\D_{\da}$}
\bl{lem-Ddown} We set for any $z>0$, $\D_{\da}(z)=\{x: l_n(x)\le z n^b\}$.
Then, for any $y_{\da}>0$, we have 
\be{lim-lemDdown}
\varlimsup_{z\to 0} \varlimsup_{n\to\infty}
\frac{1}{n^a}\log P\left(\sum_{x\in\D_{\da}(z)} 
l_n(x) \eta(x)>ny_{\da}\right)=-\infty.
\ee
\el
\bpr
We fix $z>0$ and $\{l_n(x):\ x\in\in \D_{\da}(z)\}$ 
and integrate over the $\eta$, to obtain for $\l\ge 0$
\be{equb.7}
P_{\eta}\left(\sum_{x\in\D_{\da}(z)} l_n(x) \eta(x)> n y_{\da}\right)
\le \exp(-ny_{\da}\frac{\l}{zn^{b}}) 
\prod_{x\in \D_{\da}(z)} E_{\eta}[\exp\left( \l \eta(x)\frac{l_n(x)}{zn^b}\right)].
\ee
Note that by hypothesis \reff{eq1.1}, there is $\delta>0$ such that
\be{equb.4}
\nu(\delta)=E_{\eta}[\eta^2(x) \exp(\delta |\eta(x)|)<\infty.
\ee
Also, it is an obvious fact that for $0\le \theta \le \delta$
\be{equb.5}
\exp(\theta \eta(x))\le 1+\theta \eta(x)+\frac{\theta^2 \eta(x)^2}{2}
e^{\delta |\eta(x)|}.
\ee
Thus, after taking expectation in \reff{equb.5}
\be{equb.6}
E_{\eta}[\exp(\theta \eta(x))]\le 1+\frac{\theta^2}{2} \nu(\delta)\le
e^{\theta^2 \nu(\delta)/2}.
\ee
Back to estimating \reff{equb.7}, we choose $\l\le \delta$ and use \reff{equb.6}
to obtain
\be{equb.8}
\prod_{x\in \D_{\da}(z)} E_{\eta}[\exp\left( \l \eta(x)\frac{l_n(x)}{zn^b}\right)]\le
\exp\left( \frac{\l^2}{2} \nu(\delta)
\frac{\sum_{x\in\D_{\da}(z)}l_n(x)^2}{(zn^b)^{2}}\right).
\ee
Now, $\sum_{x\in\D_{\da}(z)}l_n(x)^2\le zn^{1+b}$. Thus, 
\be{equb.9}
P\left(\sum_{x\in\D_{\da}(z)} l_n(x) \eta(x)> n y_{\da}\right)
\le\exp\left( -\frac{n^{1-b}}{z} 
\sup_{0\le \l \le \delta} \{ y_{\da} \l -\frac{\nu(\delta) \l^2}{2}\}\right)
\ee
The results follows since for any $y_{\da}>0$, the
supremum is positive, and $z$ can be sent to zero.
\epr
\subsection{Contributions of $\D_{\ua}$}
\bl{lem-Dsup} For $\D_{\ua}$ given in \ref{def-DUL}, there is $C_d>0$
such that for $n$ large 
\be{equb.21}
\P_0\left(\D_{\ua}\not= \emptyset\right)\le
C_d n e^{-\kappa_0 (yn)^a}.
\ee
\el
\bpr
First, note that 
\be{equb.19}
\P_0\left(\D_{\ua}\not= \emptyset\right)= 
\P_0\big(l_n(x)>(yn)^a\text{ for some } x\in \Ra_n\big)
\le \sum_{x\in \Z^d} \P_0(x\in \Ra_n)
\P_0\left(l_{\infty}(0)\ge(yn)^a\right).
\ee
Now, for $\kappa_0:=1/\log(\P_0(H_0<\infty))$, it is clear that
\be{equb.20}
\P_0\left(l_{\infty}(0)\ge(yn)^a\right)\le e^{-\kappa_0 (yn)^a}.
\ee
Thus, we conclude by recalling standard estimates \`a la Dvoretzky-Erd\"os (see 
for instance Theorem 6.2 of \cite{lawler}) which establish that $\E_0[\Ra_n]$
is of order $n$ in dimension larger than 2.
\epr
\subsection{Contributions of $\D_i$ for $i=0,\dots,N$.}\label{small}
\bl{lem-Di} Fix $i=0,\dots,N$. We have a constant $c_2>0$ such that
for $y>0$, and $n$ large
\[
P\big(\A, \sum_{x\in\D_i}l_n(x) \eta(x)> n y_i\big)\le 
\exp\left(-c_2 (ny)^a\right).
\]
\el
\bpr
We first treat the case $i>0$.
Note that $|\D_i|\le y^{-a} n^{1-b_i}$, and on $\A$ there are at most
$\binom{n^d}{|\D_i|}$ possible choices for $\D_i$ since the walk does
not exit a region of radius $n$. Thus, using Lemma~\ref{lem-handy} (and
\reff{ineq-use})
\ba{ineq-base}
P\big(\A,\sum_{x\in\D_i}\!\!\! &\!\!\!&\!\! l_n(x) \eta(x)> n y_i\big)\le
\sum_{L=1}^{y^{-a}n^{1-b_i}} \P_0\left(\A, |\D_i|=L\right)\ 
P_{\eta}\!\!  
\left( \sum_{j=1}^L \eta_j > \frac{ny_i}{\max\{l_n(x):x\in \D_i\}}\right)\cr
&\le & \sum_{L=1}^{y^{-a}n^{1-b_i}} \binom{n^d}{L} \sup_{\Lambda: |\Lambda|=L}
\P_0(\D_i=\Lambda) 
\ P_{\eta}\!\!\left( \sum_{j=1}^L \eta_j > \frac{ny_i}{n^{b_{i+1}}y^a}\right)\cr
&\le & \sum_{L=1}^{y^{-a}n^{1-b_i}} (n^d)^{L} \sup_{\Lambda: |\Lambda|=L}
\P_0(l_{\infty}(\Lambda)> L y^an^{b_i})
\ P_{\eta}\!\!\left(\sum_{j=1}^L \eta_j > \frac{ny_i}{n^{b_{i+1}}y^a}\right).
\ea
By using Lemma~\ref{lem-amin}, we have
\be{equb.11}
(n^d)^{L} \sup_{\Lambda: |\Lambda|=L}
\P_0\left(l_{\infty}(\Lambda)> L y^an^{b_i}\right)\le \exp\left( -\kappa_d y^a
n^{b_i}L^{1-\frac{2}{d}}+L\log(n^d)\right),
\ee
and the combinatorial factor $n^{dL}$ is negligible when
\be{equb.12}
n^{b_i}L^{1-\frac{2}{d}}\gg L\log(n).
\ee
Since $L\le y^{-a}n^{1-b_i}$ and $b_i\ge b=1/(\a+1)$, \reff{equb.12} requires $n$
large and 
\be{equb.13}
\frac{2}{d}(1-b_i)< b_i\Longleftarrow
\a< \frac{d}{2}.
\ee
Thus, the combinatorial factor is always innocuous when $\a<d/2$.

Let $B>0$ be a fixed large constant.
We say that $L$ is {\it large} when $n^{b_i} L^{1-\frac{2}{d}} > B n^{a}$,
and this case poses obviously no problem since the
term $\P_0(l_{\infty}(\Lambda)> L n^{b_i}y^a)$ suffices
to obtain the right speed. Thus, we assume that $L$ is {\it small},
that is:
\be{equb.14}
n^{b_i} L^{1-\frac{2}{d}}\le B n^{a}.
\ee
Thus, we consider for a fixed $i=1,\dots,N$
\be{equb.15}
L\le A n^{\gamma_i},\quad\text{with}\quad
\gamma_i:= \frac{a-b_i}{1-\frac{2}{d}},\text{ and }
A=B^{\frac{1}{1-2/d}}
\ee
We want to evaluate $P_{\eta}( \sum_{\D_i} \eta(x) > n^{1-b_{i+1}} y_iy^{-a},
|\D_i|\le L)$
when $L$ is as in \reff{equb.15}. First, note that
$n^{\gamma_i}\ll n^{1-b_{i+1}} y_i$,
when $n$ is large enough. Indeed, first rewrite
\[
1-b_{i+1}-\gamma_i=\left(1-a-\frac{\gamma_i}{\alpha}\right)+
\left(a-b_{i+1}+(1-\delta_0)(a-b_i)\right).
\]
Then, by noting that $1-a-\gamma_i/\a\ge 1-a-\gamma_1/\a\ge b \d_0$,
and using \refeq{ineq.req1},
\be{ineq-enough}
n^{1-b_{i+1}-\gamma_i} y_i = n^{1-a-\frac{\gamma_i}{\alpha}} 
n^{a-b_{i+1}+(1-\delta_0)(a-b_i)} y_i
\geq  n^{1-a-\frac{\gamma_1}{\alpha}} \beta y =n^{b \delta_0} \beta y .
\ee
Hence for $L$ satisfying \refeq{equb.15}, $n^{1-b_{i+1}} y_i \gg L$, and
using standard Large Deviations
estimates (see lemma \ref{TGD.lem} in the appendix),
for all $\epsilon > 0$ and
$n$ sufficiently large,
\ba{equb.16}
P_{\eta} \left( \sum_{j=1}^L \eta_j > n^{1-b_{i+1}} \frac{y_i}{y^a}\right)
&\le& \exp\left( -c_{\a}(1-\epsilon)
L\left(\frac{n^{1-b_{i+1}}y_i}{Ly^a}\right)^{\a}\right)\cr
&\le& \exp \pare{-c_{\alpha}(1-\epsilon)A
     \pare{\frac{n^{1-b_{i+1}-\gamma_i(1-1/\alpha)} y_i} {A y^a}}^\alpha}\cr
&\le& \exp \pare{-c_{\alpha}(1-\epsilon)
\frac{\beta^{\alpha}}{A^{\alpha-1}} (ny)^a} \, ,
\mbox{ using \refeq{ineq.req1}.}
\ea

By the same arguments, we can treat the case $\D_0$. Indeed, note that
$\gamma_0=(a-b)(1-2/d)^{-1}<a$, and for $L$ {\it small}, we have
\be{equb.18}
P_{\eta} \left( \sum_{j=1}^L \eta_j > n^{1-b} \frac{y_0}{y^a}\right)\le
\exp\left( -c_{\a}\frac{n^{a\a}(y_0/y^a)^{\a}}{L^{\a-1}}\right),
\ee
which is negligible since $a\a-(\a-1)\gamma_0>a$.
\epr

\br{rem-a1} When $\a=1$, then $b=a=1/2$. Thus, the range of the walk 
is divided into three sets: $\D_{\da},\D_{\ua}$ and $\D_0$ as in
\reff{def-DUL} and \reff{def-D0} respectively. Also, we can choose
$y_{\da}=y_{\ua}=y_0=y/3$. Now, our treatment
for $\D_{\da},\D_{\ua}$ only assumed small exponential moments for the walk,
which hold in this case.
To treat $\D_0$ note that only the case $L$ {\it small} may pose problem.
However, since $\gamma_0=0$, $L$ small means $L\le A$ for a large
constant. It is easy to see that those terms are of the correct order
since, there is a constant $\bar c$ such that
\be{equb.a1}
P_{\eta}\left(\sum_{j=1}^L \eta_j> n^{1-a}\frac{y_0}{y^a}\right)\le
LP_{\eta}(\eta_1>\frac{(ny)^{1/2}}{3L})\le L\exp(-\bar c (ny)^{1/2}).
\ee
\er
\section{Appendix}
\subsection{On bell-shaped densities}
We recall that a density $f$ is bell-shaped if it is
even and decreasing on $\R^+$. Our first observation is
the following.
\bl{lem-conv} If $f,g$ are two bell-shaped densities, so is
their convolution $f*g$.
\el
\bpr
First, it is obvious that $f*g$ is even. Indeed, by the eveness of
both $f$ and $g$
\[
f*g(-t)=\int_{\R}\!\! f(-t-s)g(s)ds=\int_{\R}\!\! f(t+s) g(-s) ds=
\int_{\R}\!\! f(t-s)g(s)ds= f*g(t).
\]
Now, assume that $f$ is differentiable. Then,
\be{eq2.9}
(f*g)'(t)=\int_{\R}\!\! f'(t-s)g(s)ds=\int_{\R}\!\! f'(s) g(t-s)=
\int_0^{\infty}\!\! f'(s)\left(g(|t-s|)-g(t+s)\right) ds,
\ee
where we used the oddness of $f'$ and the eveness of $g$.
Now, for $t,s\ge 0$, we have $g(|t-s|)-g(t+s)\ge 0$, and $f'(s)\le 0$
implying that $(f*g)'(t)\le 0$.

Now let $\{\v_{\e},\e>0\}$ be a differentiable {\it bell-shaped} approximate identity.
By what we just saw, $\v_{\e}*f$ is a bell-shaped differentiable density. 
So is in turn $(\v_{\e}*f)*g$. Thus, for any $0\le t\le T$, we have
$(\v_{\e}*f)*g(t)\ge (\v_{\e}*f)*g(T)$. By pointwise convergence, 
as $\e$ tends to 0, we obtain that $f*g(t)\ge f*g(T)$.
\epr

By induction, using Lemma~\ref{lem-handy}, we obtain the following corollary.
\bc{cor-obvious} If $\{\eta_i,i=1,\dots,n\}$ are independent bell-shaped
variables and $S=\a_1\eta_1+\dots+\a_n\eta_n$, with positive $\{\a_i\}$,
then $S$ has a bell-shaped density.
\ec

Finally, the useful result is the following.
\bl{lem-useful} Let $\{\eta_i,i=1,\dots,n\}$ be independent bell-shaped
variables and $0\le \a_i\le \beta_i$ for $i=1,\dots,n$. then
for any $y>0$, we have \reff{ineq-handy}.
\el
\bpr
We prove the Lemma by induction on the number of $\beta_i$ larger
than $\a_i$. Thus, it is enough to show that $x\mapsto P(S+x\eta_n)$ is
increasing on $\R^+$ when $S$ is a bell-shaped variable independent of $\eta_n$.

First note that for symmetric independent $\xi,\eta$, we have for $y>0$
\be{eq2.10}
P(\xi+\eta>y)=P(\xi>y)+\int_0^{\infty} P(\eta>z)(f_{\xi}(|y-z|)-f_{\xi}(y+z))dz.
\ee
The proof is concluded as we apply \reff{eq2.10} to $\xi=S$ and $\eta= x\eta_n$,
and as we note that $f_{S}(|y-z|)-f_{S}(y+z)\ge 0$ and $x\mapsto P(\eta_n>z/x)$
is increasing.
\epr

\subsection{On a localization result}
We first prove Lemma~\ref{lem-amin}.

\noindent First Step: 
We show that $\E_0[l_{\infty}(\Lambda)]\le C_d |\Lambda|^{2/d}$.

The following Green function estimates is standard (see for instance
\cite{lawler} Theorem 10.1): there is $C_d$ such that
for any $y\in \Z^d$
\be{eq2.1}
G(0,y)=\E_0[l_{\infty}(y)]\le \frac{C_d}{1+||y||^{d-2}}.
\ee
Now, $l_{\infty}(\Lambda)=\sum_{y\in \Lambda}l_{\infty}(y)$ and
\be{eq2.2}
\E_0[l_{\infty}(\Lambda)]=\sum_{y\in \Lambda} \E_0[l_{\infty}(y)]\le
\sum_{y\in \Lambda} \frac{C_d}{1+||y||^{d-2}}.
\ee
We establish now an upper bound on the right hand side of \reff{eq2.2}.
Let $\v$ be an ordering of the sites of $\Z^d$ in increasing distance from
the origin. In other words, $\v:\N\to\Z^d$ is a one to one, onto map so that
$||\v(i)||\le ||\v(i+1)||$, for all $i\in \N$. Let $\psi:\{0,\dots,|\Lambda|-1\}
\to \v^{-1}(\Lambda)$ be an ordering of $\v^{-1}(\Lambda)$ (so that
$\psi(0)<\dots<\psi(|\Lambda|-1)$) and note that 
\be{eq2.3}
k\le \psi(k),\quad\text{and}\quad ||\v(k)||\le ||\v(\psi(k))||.
\ee
Thus, $g:=\v\psi^{-1}\v^{-1}:\Lambda\to \Z^d$ is a rearrangement of
$\Lambda$ inside a ``ball'' of radius proportional to $|\Lambda|^{1/d}$.
Thus, it is a trivial fact that there is a constant $c_d'$ and $\sup_{\Lambda}
||g(x)||\le c_d'|\Lambda|^{1/d}$. Let $r:=c_d'|\Lambda|^{1/d}$, and note that
\ba{eq2.4}
\sum_{x\in \Lambda} \frac{1}{1+||x||^{d-2}}&\le &
\sum_{x\in \Lambda} \frac{1}{1+||g(x)||^{d-2}}\le
\sum_{||y||\le r} \frac{1}{1+||y||^{d-2}}\cr
&\le & \frac{1}{2}+\int_0^r \frac{s^{d-1}}{1+s^{d-2}}ds\le
\frac{1}{2}+\int_0^r sds\cr
&\le & \frac{1}{2}+\frac{r^2}{2}\le r^2.
\ea
The first step concludes easily.
By Chebychev's inequality we have
\be{eq2.5}
\sup_{x,\Lambda} \P_x(l_{\infty}(\Lambda)>2C_d r^2)<\frac{1}{2}.
\ee
Indeed, the starting point of the walk
can very well be any site $x\in \Z^d$ since the transition kernel
is translation invariant, and $\Lambda$ is arbitrary. 

Second Step: We show that 
\[
\P_0(l_{\infty}(\Lambda)>t)\le (\frac{1}{2})^{t/(2C_d r^2)}
\]
Define a sequence of stopping times for $k=1,2,\dots$
\be{eq2.6}
\sigma_k=\inf\{n\ge 0: l_n(\Lambda)>2kC_d r^2\},
\ee
and note that $\sigma_k=\sigma_{k-1}+\sigma_1\circ \theta_{\sigma_{k-1}}$.
We have used the notation $\theta_k$ for the time translation by $k$-units.
Now, the bound \reff{eq2.5} can be expressed in term of $\sigma_1$ as
\[
\P_0(\sigma_1<\infty)<\frac{1}{2}.
\]
We express now the total sojourn time in $\Lambda$ in terms of $\{\sigma_k,k\in \N\}$
\be{eq2.7}
\P_0(l_{\infty}(\Lambda)>2kC_d r^2)=\P_0(\sigma_k<\infty),
\ee
and by the Strong Markov property
\ba{eq2.8}
\P_0(l_{\infty}(\Lambda)>2kC_d r^2)&=&\E_0[\ind_{\sigma_{k-1}<\infty}
\P_0(\sigma_{k-1}+\sigma_1\circ \theta_{\sigma_{k-1}}<\infty|\FF_{\sigma_{k-1}})]\cr
&=&\E_0[\ind_{\sigma_{k-1}<\infty}\P_{S_{\sigma_{k-1}}}(\sigma_1<\infty)]\cr
&\le & \frac{1}{2} \P_0(\sigma_{k-1}<\infty)
\ea
By induction the bound \reff{key-amin} follows readily.

\subsection{On a Large Deviation estimate.}
To be self-contained, we give an obvious estimate, for
which a reference could not be found. We assume that $\a>1$.
\bl{TGD.lem}
For all $\epsilon > 0$, $L$ a positive integer, and for $t/L$ large enough,
\[
P_{\eta} \cro{ \sum_{j=1}^L \eta_j \geq t}
 \leq \exp \pare{ -c_{\alpha} (1-\epsilon) \frac{t^{\alpha}}{L^{\alpha-1}}}
\, .
\]
\el

\bpr

For $\lambda  \in \R$, set $\Lambda(\lambda) := \log E_{\eta}
\cro{e^{\lambda \eta_1}}$.
\[
P_{\eta} \cro{ \sum_{j=1}^L \eta_j \geq t}
\leq \exp\pare{-\alpha c_{\alpha} \pare{\frac{t}{L}}^{\alpha-1} t}
 \exp \pare{L \Lambda(\alpha  c_{\alpha} (t/L)^{\alpha-1})} \, .
\]
By Kasahara's Tauberian theorem, for large $x$
\[
\Lambda(x) \simeq \frac{1}{\bar{\alpha}
(\alpha c_{\alpha})^{\bar{\alpha}-1}} x^{\bar{\alpha}},\quad\text{where}\quad 
\frac{1}{\alpha} +
\frac{1}{\bar{\alpha}} =1.
\]
Hence for all $\epsilon >0$ and $t/L$ large enough,
\beqn
P_{\eta} \cro{ \sum_{j=1}^L \eta_j \geq t}
& \le& \exp\pare{-\alpha c_{\alpha} \frac{t^{\alpha}}{L^{\alpha-1}} }
\exp \pare{(1+\epsilon) \frac{\alpha c_{\alpha}}{\bar{\alpha}}
\frac{t^{\alpha}}{L^{\alpha -1}}}\cr
& \le& \exp\pare{- c_{\alpha} \frac{t^{\alpha}}{L^{\alpha-1}}
\pare{1-\frac{\epsilon \alpha}{\bar{\alpha}}}}  \, .
\eeqn
\epr


\begin{thebibliography}{99}
\bibitem{AC1}  Asselah, A.; Castell, F.
{\it Quenched large deviations for diffusions in a random
Gaussian shear flow drift.}
Stochastic Process. Appl. {\bf 103} (2003), no 1, 1-29.

\bibitem{AC2} Asselah, A; Castell, F.
{\em Large deviations for Brownian motion in a random scenery.}
Probab. Theory Related Fields  126  (2003),  no. 4, 497--527.

\bibitem{B} Bolthausen, E.
{\em A central limit theorem for two-dimensional random walk in random
sceneries.} Ann. Probab. {\bf 17} (1989), no. 1, 108-115.


\bibitem{B1} Borodin, A. N.
{\em Limit theorems for sums of independent random variables defined on
a transient random walk.}
Investigations in the theory of probability distributions, IV.
Zap. Nauchn. Sem. Leningrad. Otdel. Mat. Inst. Steklov. (LOMI) {\bf 85}
(1979), 17-29, 237, 244.

\bibitem{B2} Borodin, A. N.
{\em A limit theorem for sums of independent random variables defined on
a recurrent  random walk.}
Dokl. Akad. Nauk. SSSR {\bf 246} (1979), no. 4, 786-787.

\bibitem{CP} Castell, F.; Pradeilles, F.
{\em Annealed large deviations for diffusions in a random
Gaussian shear flow drift.}
Stochastic Process. Appl. {\bf 94} (2001), 171--197.

\bibitem{C} Castell, F.
{\em Moderate deviations for diffusions in a random Gaussian shear
flow drift.}
Ann. Inst. H. Poincar\'e Probab. Statist. {\bf 40} (2004), no 3,
337-366

\bibitem{CFRRS} Cs\'aki, E.; F\"oldes, A.; R\'ev\'esz, P.;  Rosen, J.;
Shi, Z.
{\em Frequently visited sets for random walks.} Preprint 2004.
arXiv:math.PR/0412018

\bibitem{KS} Kesten, H.; Spitzer, F.
{\em  A limit theorem related to a new class of self-similar processes.}
Z. Wahrsch. Verw. Gebiete {\bf 50} (1979), no. 1, 5--25.

\bibitem{GKS} Gantert, N.; K\"onig, W.; Shi, Z.
{\em Annealed deviations of random walk in random scenery}
Preprint 2004. arXiv:math.PR/0408327

\bibitem{HGK} van der Hofstad, R., Gantert, N.; K\"onig, W.
{\em Deviations of a random walk in a random scenery with
stretched exponential tails} preprint 2004.\\
arXiv:math.PR/0411361

\bibitem{lawler} Lawler, G.{\em Notes on random walks}, In preparation.\\
www.math.cornell.edu/$\sim$lawler/m778s04.html

\bibitem{Nagaev} Nagaev, A. V.
{\em A property of sums of independent random variables.}
Teor. Verojatnost. i Primenen. 22 (1977), no. 2, 335--346.

\end{thebibliography}
\end{document}